\documentclass[a4paper,12pt]{article}

\newtheorem{thm}{Theorem}

\newtheorem{prop}{Proposition}

\newtheorem{cor}{Corollary}
\newtheorem{lem}{Lemma}
\newcommand{\tP}{\tilde{P}}
\newcommand{\tZ}{\tilde{Z}}
\newcommand{\tL}{\tilde{\cal L}}
\newcommand{\Hilb}{{\rm Hilb}}
\newcommand{\cK}{{\cal K}}

\newcommand{\dbd}{\overline{\partial}\partial}
\newcommand{\cM}{{\cal M}}
\newcommand{\bC}{{\bf C}}
\newcommand{\bR}{{\bf R}}

\newcommand{\cL}{{\cal L}}

\newcommand{\cG}{{\cal G}}

\title{Scalar curvature and projective embeddings, II}
\author{S. K. Donaldson}
\begin{document}
\maketitle

\section{Introduction}
This  is a sequel to the  previous paper \cite{kn:don}, which studied
connections between the  differential geometry of complex projective varieties
and certain specific \lq\lq balanced'' embeddings in projective space.
  The original
plan was that this sequel would be a  lengthy paper, discussing various
extensions and ramifications of the ideas sudied in \cite{kn:don}. However
this plan has been modified in the light of subsequent developments. On
the one hand, Mabuchi \cite{kn:M1}, \cite{kn:M2}, \cite{kn:M3} has 
extended the results of \cite{kn:don} to the
case where the varieties have infinitesimal automorphims. On the other
hand, Phong and Sturm \cite{kn:PS1}, \cite{kn:PS2} have sharpened some of the arguments in
\cite{kn:don}. They also  explain the relation of the ideas to the 
 Deligne pairing and
the  Chow norm, and to earlier work of Zhang \cite{kn:Z}, which the author
was unfortunately not aware of when writing \cite{kn:don}. These developments
mean that some of the results planned for the sequel are now redundant,
while on the other hand the exposition of  all the different points of view
 has grown into a daunting task. Thus, instead, this sequel is a short
 paper  devoted  to the proof of one result which is quite an easy
 consequence of the main theorem in \cite{kn:don}.

   To state our result, suppose that $X$ 
   is a compact Kahler manifold and fix a Kahler class $[\omega_{0}]$ on $X$.
    Recall that the {\it Mabuchi functional} \cite{kn:M}
   is a functional $\cM$, defined up to an arbitrary additive constant, on the Kahler metrics
   in this cohomology class which is characterised by the formula
   \begin{equation}  \delta \cM = \int_{X}(S - \hat{S}) \ \delta \phi \frac{\omega^{n}}{n!}
   \label{eq:Mdef}\end{equation}
   Here metrics $\omega$ in the class are represented by Kahler potentials
   $\omega= \omega_{0}+ i\dbd \phi $ and $\delta \phi$ is an infinitesimal
   variation in $\phi$. The symbol  $S$ denotes
    the scalar curvature of the metric $\omega$ and $\hat{S}$ is the average
    value of $S$ with respect to the volume form $\frac{\omega^{n}}{n!}$,
    which is a topological invariant of the Kahler class. What equation
    (1)
    really defines is a $1$-form on the space of metrics in the Kahler
    class and one shows that this is closed, so is the derivative of a
    function $\cM$, unique up to a constant. 
   
    Now suppose that $L$ is a positive line bundle over $X$ and the Kahler
    class is $2\pi c_{1}(L)$. As in \cite{kn:don},
    we write ${\rm Aut}(X,L)$ for the group of automorphims of the pair
    $(X,L)$ modulo the trivial automorphisms $\bC^{*}$ (acting by constant
    scalar multiplication on the fibres). 
    
     \begin{thm}
       Suppose that ${\rm Aut}(X,L)$ is discrete and that there is a metric
      $\omega_{\infty}$  of constant scalar curvature in the Kahler class
      $2\pi c_{1}(L)$. Then $\omega_{\infty}$ minimises the Mabuchi functional
      in this Kahler class.
      \end{thm}
      
      No doubt the hypothesis on the automorphism group can be relaxed,
      using the techniques of \cite{kn:M1}. In the case of Kahler-Einstein
      metrics this result was proved by Bando and Mabuchi in \cite{kn:BM}.

       Before beginning the proof we include some general discussion. It
       is clear from the definition that a constant scalar curvature metric
       $\omega_{\infty}$
       is a critical point of $\cM$ and it is also easy to see that the
      it is a local minimiser. Thus the difficulty in proving the theorem
       is to go from this local picture to the space of all Kahler metrics.
       In fact it is more convenient to work with the space
       $\cK$ of Kahler potentials
     The main theme of the proof is the notion of {\it convexity}, for
     the functional $\cM$ and various related functionals. This is related to the  uniqueness of the constant
            scalar curvature metric which was the main result of \cite{kn:don}.
            The notion of convexity depends on a suitable geometric
            structure on the space $\cK$. One obvious structure is realised
            by the embedding of $\cK$ as an open set in the vector space
            $C^{\infty}(X)$, and convexity has a meaning in this sense. However,
            as far as the author knows, there is no reason to think that
            $\cM$ is convex in this obvious sense. A more subtle geometric
            structure on $\cK$ arises from considerations of symplectic
            geometry, yielding a Riemannian metric making $\cK$ formally
            into a symmetric space,  dual to a group $\cG$ of symplectomorphisms
            of $X$ (or, more precisely, the extension of the group of exact
            symplectomorphisms by $S^{1}$) (\cite{kn:M0}, \cite{kn:Semmes},
            \cite{kn:D0}). Convexity has a meaning in
            this sense, as convexity along the geodesics in $\cK$, and
            indeed $\cM$ is convex from this point of view. Moreover this
            convexity property arises as part of a general package of results,
            related to the fact that the scalar curvature furnishes a moment
            map for an action of the symplectomorphism group. As explained in
             \cite{kn:D0}, if one knew that any two points in $\cK$
            can be joined  by a geodesic segment then Theorem 1 (without
            any restrictive hypotheses) would be
            an immediate consequence of this convexity. The 
             construction of geodesic segments in $\cK$ is a very difficult
             analytical problem but great progress in developing this
              approach
             has been made by X.X.Chen. While it is perhaps not  likely that all
             points can in fact be joined by {\it smooth} geodesics, in
             \cite{kn:Chen} Chen proved an existence theorem for $C^{1,1}$
             geodesics which enabled him to deduce the bound on the Mabuchi
             functional in the case when $\hat{S}$ is negative. Chen's approach 
             has been extended
             in recent work of Chen and Tian \cite{kn:CT}, who prove a more
             general result, subsuming our Theorem 1. However we hope that
               the  proof here will still
             have some interest.

             Our proof of Theorem 1 develops the approach of \cite{kn:don}.
             There are two main parts of the proof. In the first part (Section
             2) we  consider another functional on the space
             $\cK$ which can be related to a finite-dimensional problem
             on the symmetric space $M=GL(d,\bC)/U(d)$. The main result
             here is Theorem 2 below. In the second 
             part (Section 3) we replace the line bundle $L$ by $L^{k}$ for large $k$
             and apply
              the main result of \cite{kn:don}, and the asymptotic
             theory used there. In the first part of the we
               emphasise the symmetry in the discussion between
             the spaces $M$ and $\cK$. All of this can be fitted in to
             the general picture of \cite{kn:don}, which started from 
             commuting
             actions of the two groups $U(d)$ and $\cG$. With the exception
             of Lemma 4, the picture we develop in Section 2 largely follows
             from these group actions and general constructions concerning
             symplectic quotients. However in order to prove Lemma 4 we
             need to have  the relevant formulae written down  explicitly,
             and all-in-all we make a fresh start, so that this paper can
             be read 
             independently of 
              \cite{kn:don}. The proof of Lemma 4 hinges on the convexity
               of another functional $I$ on $\cK$ in the first, more obvious,
               sense.

             \section{}
             
             Let $X,L$ be as above and suppose that we have a fixed subspace
             $E\subset H^{0}(X;L)$, of dimension $d$, which gives a projective embedding of
             $X$.    Let $h$ be a Hermitian metric on the line
             bundle $L\rightarrow X$. This defines a unitary 
             connection on $L$, with curvature $-i\omega_{h}$ say, and we have
             \begin{equation} \omega_{e^{\phi} h} = \omega_{h} + i\dbd \phi.\end{equation}
             Let $\cK$ denote the set of metrics $h$ such that 
             $\omega_{h}$
             is a positive $(1,1)$ form on $X$. Thus if we fix a base point
             $h_{0}$ in $\cK$   we can identify this space with the set of
             Kahler potentials,
             \begin{equation} \cK \cong\{ \phi:
              \omega_{h_{0}}+i \dbd \phi>0\}. \end{equation}
             We write
             $$  d\mu_{h} = \frac{\omega_{h}^{n}}{n!}: $$
             the standard volume form associated to the Kahler 
             metric $\omega_{h}$, and we write $V$ for the volume of $X$,
             computed using any of these metrics.

             Now let $M$ denote the set of Hermitian metrics on the finite-dimensional
             vector space  $E$. We
             have two fundamental constructions:
             \begin{itemize}
              \item Given $h\in \cK$ we get a metric $\Hilb(h)\in M$, the
              standard $L^{2}$-metric defined using the fibre metric on
              $L$ and the volume form $ d \mu_{h} $ but with a convenient
              normalising factor
              $$  \Vert s\Vert^{2}_{\Hilb(h)} = \frac{d}{V} \int \vert s\vert_{h}^{2}
              d\mu_{h}. $$
               So we have a map
              $$  \Hilb: \cK\rightarrow M. $$ 
              \item Given $H \in M$ we define a metric $FS(H)$ on $L$
              by decreeing that if $s_{\alpha}$ is an orthonormal basis
              of $E$ with respect to $H$ then
              $$ \sum_{\alpha} \vert s_{\alpha}\vert^{2}_{FS(H)} = 1, $$
              pointwise on $X$. The corresponding form $\omega_{FS(H)}$ is
              the restriction of the standard Fubini-Study metric on ${\bf
              P}(E^{*})$ to the image of $X$ under the projective embedding,
              so it is a positive form and we have a map
              $$   FS: M \rightarrow \cK. $$ 
              \end{itemize}

             A \lq\lq balanced metric'', in the sense of \cite{kn:don},
             corresponds to a pair $(h^{*},H^{*})$ with
             \begin{equation} \Hilb(h^{*})= H^{*}\ , \ FS(H^{*})= h^{*}, \end{equation}
             that is, to a fixed point of  either composite $FS\circ
             \Hilb$ or $\Hilb\circ FS$. 
            
            Suppose we fix a non-zero element
             $\Theta\in \Lambda^{\dim E} E$. Then we can define the determinant
             of any metric $H$ in $M$. (This is just the determinant of the matrix
             of $H$ with respect to any basis $(s_{\alpha})$ of $E$ with
             $\Theta= s_{1}\wedge\dots \wedge s_{d}$).
             Thus we have a map 
             $$ \log \det: M \rightarrow \bR.$$
             A different choice of $\Theta$ just changes this map by the
             addition of a constant. We need the analogous map on $\cK$.  
              We recall that there is a functional $I: \cK \rightarrow
              \bR,$  unique up to an additive
               constant, characterised by the condition that if $h_{t}=
               e^{\phi_{t}} h_{0}$ is  a path in $\cK$ then
            \begin{equation} \frac{d I}{dt} = \int_{X} \dot{\phi}\ 
            d\mu_{h_{t}} \label{Idef}
            \end{equation}
            This functional has been used by many authors. As for the Mabuchi
            functional, essentially what is defined by equation (5) is a
               $1$-form on $\cK$, and one has to check that this is $1$-form
               is closed. Notice that the derivative of the volume form
               is given by
               \begin{equation} \left(\frac{d}{dt} \right) d\mu_{h_{t}} = 
               \Delta'_{t}
                ( \dot{\phi})\  d\mu_{h_{t}},
               \end{equation}
               where $\Delta'_{t}$ denotes {\it one half} the ordinary
               Riemannian Laplacian of the metric
                $\omega_{h_{t}}$. (We use the sign convention in which the
                Laplacian is a nonegative operator.)
               This gives 
               \begin{equation} \frac{d^{2}I}{dt^{2}} =
                \int_{X} \left(\ddot{\phi} 
               +\dot{\phi} \Delta' \dot{\phi}  \right) d\mu_{h_{t}} \end{equation}
            \begin{lem} For any $h, h_{0} \in \cK$, with $h=e^{\phi} h_{0}$
            we have
              $$ \int_{X} \phi \ d\mu_{h_{0}} \leq
               I(h) - I(h_{0}) \leq \int_{X} \phi
               \ d\mu_{h}. $$
              \end{lem}
              This is also well-known, and amounts to the fact that 
             $I$ is a {\it convex function} on $\cK$, regarded as an open
             subset of $C^{\infty}(X)$. By symmetry it suffices to prove
             the left-hand inequality. We write $h_{t} = e^{t\phi} h_{0} $
              and consider the
             function 
             $$ f(t) = 
             \int_{X} t\phi\  d\mu_{h_{0}} - (I(h_{t}) - I(h_{0})) $$
             on $[0,1]$. Both $f(0)$ and $f'(0)$ vanish and
             $$  f''(t) = -\int_{X} \phi \Delta'_{t}
              \phi d\mu_{h_{t}} \leq
             0, $$
             so $f(1)\leq 0$.                     

                 \

               Clearly we have the following scaling identities, for $\alpha\in
               \bR$,
               \begin{eqnarray} \Hilb(e^{\alpha} h) &=& e^{\alpha} \Hilb(h)\
               \ , \ \ 
               FS(e^{\alpha} H)= e^{\alpha}FS(H),\\
               \log \det (e^{\alpha} H)&=& \log \det H + \alpha d\ \  ,\
               \ 
                I(e^{\alpha} h) = I(h)
               + \alpha V . \end{eqnarray}

               We now define maps
               \begin{eqnarray} \cL&=& \log \det \circ \Hilb:
                \cK\rightarrow \bR, \\
                Z&=& -I \circ FS: M\rightarrow \bR. \end{eqnarray}
               We also put
               \begin{equation}
               \tL= \cL- \frac{d}{V} I\ , \ \tZ = Z+ \frac{V}{d} \log \det\end{equation}
               so $\tL$ and $\tZ$ are unchanged by constant scaling on
               $\cK$ and $M$ respectively.
               
               \

                Given a metric $h\in \cK$, we define a function $\rho_{h}$ on $X$
 by
 \begin{equation} \rho_{h} = \frac{d}{V} \sum_{1}^{d} \vert s_{\alpha} \vert_{h}^{2}
, \end{equation}
where $( s_{\alpha})$ is an orthonormal basis 
 for $E$ with respect to $\Hilb (h)$. The function
 $\rho_{h}$ does not depend on the choice of this orthonormal basis.

\begin{lem}
     The derivative of $\cL$ on $\cK$ is given by
     $$ \delta\cL =   \int_{X} (\Delta' \rho_{h} + \rho_{h}) (\delta \phi). $$
     \end{lem}
     To see this note first that changing the element $\Theta$ changes
     $\cL$ by a constant and does not affect the derivative of $\cL$. So to
     verify the formula at a given point $h_{0}\in \cK$ we may suppose
     that there is an orthonormal basis $(s_{\alpha})$ with respect to 
     $\Hilb(h_{0})$ such that $s_{1}\wedge\dots s_{d}= \Theta$.
     For any $\phi$ we have 
     $ \cL(e^{\phi} h_{0}) = \log \det H_{\phi} $ where $H_{\phi}$ is
      the matrix with
     entries $\langle s_{\alpha}, s_{\beta}\rangle_{\Hilb(e^{\phi} h_{0})}$. Differentiating
     at the given point $\phi=0$, where $H_{\phi}$ is the identity matrix, we have
     $$ \delta \cL = \sum_{\alpha} \delta \Vert s_{\alpha} 
     \Vert_{\Hilb(e^{\phi}h_{0})}^{2}. $$
     Now, pointwise on $X$,
     $$\delta \vert s_{\alpha}\vert_{e^{\phi}h_{0}}^{2}=
      (\delta \phi)\vert s_{\alpha}\vert_{h_{0}}^{2},
     .$$
     When we take the $L^{2}$ norm over $X$ we get another term from the
     variation in the volume form (3); thus
     $$   \delta \Vert s_{\alpha}\Vert_{\Hilb(h)}^{2} = \frac{d}{V} \int_{X} 
     (\delta \phi
     + \Delta' (\delta \phi)) \vert s_{\alpha}\vert_{h_{0}}^{2}\  d\mu_{h_{0}}.
     $$Summing over $\alpha$ and integating by parts we obtain
     \begin{eqnarray} \sum_{\alpha} \delta \Vert s_{\alpha} \Vert_{\Hilb(h)}^{2} 
     &=& \int_{X}
     \left( \delta \phi + \Delta'(\delta \phi) \right) \rho_{h_{0}}
      d\mu_{h_{0}}\nonumber\\ &=& 
     \int_{X} \left( \rho_{h_{0}} + \Delta' \rho_{h_{0}}\right) \delta
     \phi\  d\mu_{h_{0}}, \nonumber\end{eqnarray}
     as required.

     \begin{cor}  A  point $h^{*}$ in $\cK$ is balanced if and only if
     it is a  
      critical point of the functional 
      $\tL= \cL- \frac{d}{V}I$ on $\cK$.
       \end{cor} This is just because $\delta
      \cL$ vanishes for all $\delta \phi$ of integral zero if and only
      if $\Delta' \rho_{h}+ \rho_{h}$ is a constant, but this can only occur
      if $\rho_{h}$ is a constant, since $\Delta'$ is a non-negative operator.
      The factor inserted in the definition of $\Hilb$ then implies that
      $h=FS(Hilb(h))$. 
      
      Symmetrically we have
      \begin{lem} The derivative of $Z$ at a point $H$ of  $M$ is given by
      $$ \delta Z = \sum_{\alpha\ \beta} \int_{X}  (\delta H)_{\alpha \beta}
      (s_{\alpha}, s_{\beta})_{FS(H)} d\mu_{FS(H)}, $$
      where $(s_{\alpha})$ is an orthonormal basis of $E$ with respect to
      $H$. \end{lem}
      The proof is completely straightforward. It follows then that a balanced
      point (regarded as fixed point of $\Hilb \circ FS$ on $M$)
      is the same as a critical point of the function $\tZ= Z+ \frac{V}{d} \log
      \det$ on $M$.

      We can now state the main result of this section,
      \begin{thm}
      Suppose that there is a  balanced point $h$ in $\cK$. Then $h$
      is an absolute minimum of the functional $\tL$
      on $\cK$.
      \end{thm}

        To prove this we
        consider the function $P$ on the product $\cK \times M$ given by
       \begin{equation}  P(h,H) = \log\ {\rm Tr}\ ( \Hilb(h) H^{-1}). \end{equation}
       In other words 
       $$  P(h,H) = \log \sum_{\alpha} \Vert s_{\alpha}\Vert^{2}_{\Hilb(h)},
       $$
       where $(s_{\alpha})$ is an $H$-orthonormal basis of $E$. From the
       definitions, we
       have
       $$  P(h, \Hilb(h))=  \log d\ \ , \  P(FS(H), H)= \log d, $$
       for all $h\in \cK, H\in M$.
      
      We put
      $$ \tP(h,H) = P(h,H) -\log d + \frac{1}{d} \log \det H - \frac{1}{V}
      I(h), $$
      then 
      \begin{equation}\tP(h, \Hilb(h)) = \frac{1}{d} \tL(h)\ ,\  \tP(FS(H), H)= \frac{1}{V}
      \tZ(H). \end{equation}

       We now turn to the finite-dimensional side, with the function
       $Z$ on $M$, the set of Hermitian metrics on $E$. There is a standard
       metric on $M$, regarded as a symmetric space $GL(d,\bC)/U(d)$
        and  the geodesics in this metric are just the images of the
       $1$-parameter subgroups in $GL(d,\bC)$.  The crucial result we need
       is 
       \begin{prop}
       The function $Z$ is convex along geodesics in $M$.
       \end{prop}
       This result is essentially equivalent to that of Zhang 
       \cite{kn:Z}, and is 
       also implicit in  \cite{kn:don}. In the framework of
       \cite{kn:don} this is an instance of the convex function one always obtains associated
       to a Hamiltonian action of a group on a Kahler manifold (see \cite{kn:DK},
       Chapter 6, for example). This is the well-known \lq\lq Kempf-Ness
       principle'' \cite{kn:KN}. However we give a direct proof here since
       it  is quite short (once one knows the relevant identities to use).
       
       A general geodesic in $M$ can be written as
       $$ H_{t} = {\rm diag}\ ( e^{\lambda_{\alpha} t} ), $$
       with respect to some  basis $s_{\alpha}$, orthonormal with respect
       to $H_{0}$. By definition
       $$ Z(H_{t}) =- I( h_{t}), $$
       where $h_{t} = FS(H_{t})$. Thus applying (7), and evaluating at $t=0$,
       $$ \frac{d^{2} Z}{dt^{2}} = -\int_{X}\dot{\phi} \Delta' \dot{\phi}
        + \ddot{\phi}\ 
       d\mu_{h_{0}}. $$
      
       Now  the sections 
       $$ \exp(-\frac{\lambda_{\alpha}}{2} t) s_{\alpha}$$ are
       orthonormal with respect to $H_{t}$ so by definition $h_{t}=e^{\phi_{t}}
       h_{0}$ where
       $$ \phi_{t} =
        -\log \sum_{\alpha} e^{-\lambda_{\alpha} t} \vert s_{\alpha}\vert_{h_{0}}^{2}.
        $$
       Then, using the fact that $\sum_{\alpha} \vert s_{\alpha}
       \vert_{h_{0}}^{2}= 1$ we have, at $t=0$,
       $$ \dot{\phi} = \sum \lambda_{\alpha} \vert s_{\alpha}\vert^{2},
       $$
       $$  \ddot{\phi} = -\sum \lambda_{\alpha}^{2} \vert s_{\alpha}\vert^{2}
       + \left( \sum \lambda_{\alpha} \vert s_{\alpha}\vert^{2}\right)^{2}.
       $$
       
       We have then
       $$ \ddot{Z}=\int_{X} -\frac{1}{2} \vert \nabla \dot{\phi}\vert^{2} 
       + \sum_{\alpha} \lambda_{\alpha}^{2} \vert s_{\alpha}\vert_{h_{0}}^{2} -
       \left( \sum_{\alpha} \lambda_{\alpha} \vert s_{\alpha}\vert_{h_{0}}^{2}\right)^{2}.
       $$
       We  denote by $(\ ,\ )$ any of the four natural $\bR$-bilinear pairings
       \begin{eqnarray} T^{*}X&\times& (T^{*}X\otimes L)\rightarrow L\ \
       \ \ \ \ 
        L\times
       (T^{*}X \otimes L) \rightarrow T^{*}X \nonumber\\
       L&\times& L \rightarrow \bR\ \ \ \ \ \ \ \ \ \ \ \ \ \ \ \ \ \  \ T^{*}X\times T^{*}X \rightarrow
       \bR\nonumber\end{eqnarray} 
       obtained using the metrics on $L$ and $T^{*}X$ furnished by $h_{0}$.
       The crucial identity we need is that for any function $f$ on $X$,
       \begin{equation} \vert \nabla f \vert^{2}= 2 \sum_{\alpha} \vert (
       \nabla f , \nabla s_{\alpha} )\vert_{h_{0}}^{2}. \end{equation}
       This follows, with a little thought, from the definition of the
       metric $FS(H_{0})$---the restriction of the Fubini-Study metric,
       see \cite{kn:don} page 502, but note that there is an error there
       in the factor of $2$ above. Consider the positive function
       $$ F= \sum_{\alpha} \vert (\nabla \dot{\phi}, \nabla s_{\alpha})
       - (\lambda_{\alpha} - \dot{\phi}) s_{\alpha}\vert_{h_{0}}^{2} .
       $$on $X$. Expanding out, and applying (16) with $f=\dot{\phi}$,
       $$F=  \frac{1}{2} \vert \nabla \dot{\phi}\vert^{2}
       + \sum_{\alpha}(\lambda_{\alpha} - \dot{\phi})^{2}\vert 
       s_{\alpha}\vert_{h_{0}}^{2}
       -2 ((\nabla\dot{\phi}, 
       \nabla s_{\alpha}), s_{\alpha}) (\lambda_{\alpha} - \dot{\phi}). $$
       Now  we have an identity
       $$ ((\nabla\dot{\phi}, \nabla s_{\alpha}), s_{\alpha}) = (\nabla
       \dot{\phi}, (s_{\alpha}, \nabla s_{\alpha}))
       $$ and $\nabla \vert s_{\alpha}\vert_{h_{0}}^{2}= 2 (s_{\alpha}, \nabla s_{\alpha})$
       so
       $$2((\nabla \dot{\phi},\nabla s_{\alpha}), s_{\alpha}) (\lambda_{\alpha}
       -\dot{\phi}) =  (\nabla \dot{\phi}, \nabla \vert
        s_{\alpha}\vert_{h_{0}}^{2})(\lambda_{\alpha}-
       \dot{\phi}). $$
       Now when we sum over $\alpha$ we can use the relations
       $$ \sum \vert s_{\alpha}\vert^{2}=1\ , \sum \lambda_{\alpha} \vert
       s_{\alpha}\vert^{2} = \dot{\phi}, $$ to obtain
       
       $$ F= \frac{1}{2} \vert \nabla \dot{\phi} \vert^{2} 
        + \sum (\lambda_{\alpha} - \dot{\phi})^{2} \vert
       s_{\alpha}\vert_{h_{0}}^{2}- \vert \nabla \dot{\phi} \vert^{2}. $$
 
       and so  finally we see that
       $$\ddot{Z}= \int_{X} F\  d\mu_{h_{0}}$$
       (using again the condition $\sum \vert s_{\alpha}\vert^{2}=1$).
       Thus $\ddot{Z}\geq 0$ as required.

       \begin{cor} If $H^{*}\in M$ is a balanced point then it is an absolute
       minimum of the function $\tZ=Z-\frac{V}{d} \log \det$ on $M$.
       \end{cor}
       To see this note that that $\log \det $ is linear on geodesics,
       so $Z-\frac{V}{d} \log \det$ is also convex on geodesics in $M$.
       Since any two points in $M$ can be joined by a geodesic any critical
       point is an absolute minumum.
       
       (Just as Proposition 1 can be obtained from the group actions considered
       in \cite{kn:don} and general theory, so also one can see that the
       symmetrical result holds: the functional $\cL$ is convex along geodesics
       in $\cK$. However we will not go into this since the whole point
       of our proof is to avoid using the geodesics in $\cK$.)

       The other ingredient in the proof of Theorem 2 is the next Lemma.
       
       \begin{lem} For any $h\in \cK$ and $H\in M$:
          $$ \tP(h, H)\geq \tP(FS(H), H). $$
             \end{lem}

    Choose an orthonormal basis $s_{\alpha}$ with respect to $H$ and
    let $h_{0}= FS(H)$, so $\sum \vert s_{\alpha}\vert_{h_{0}}^{2}=1$.
    Let $h=e^{\phi} h_{0}$ so
    $\sum \vert s_{\alpha}\vert_{h}^{2} = e^{\phi}$ and
    $$ \sum \Vert s_{\alpha}\Vert^{2}_{\Hilb(h)} =\frac{d}{V}\int_{X} e^{\phi} d\mu_{h}.
    $$Then, going back to the definitions, we find 
    $$
     \tP(h,H)- \tP(FS(H),H)= \log \left(\frac{1}{V} \int_{X} e^{\phi} d\mu_{h}\right)
    - \frac{1}{V} ( I(h)- I(h_{0}))$$
                        
                         Now, by the convexity of the exponential function
                   on $\bR$,
                   $$  \log \left( \frac{1}{V} \int_{X} e^{\phi}
                   d\mu_{h}\right) \geq \frac{1}{V} \int_{X} \phi\ d\mu_{h},
                   $$
                   and by Lemma 1
                   $$ \int_{X} \phi d\mu_{h} \geq I(h) - I(h_{0})$$
                   which completes the proof.

    The proof of Theorem 2 is now in our hands. We have
   
    \begin{equation} \tP(FS(H), H)= \log d + \frac{1}{V} \tZ. \end{equation}
    Suppose $(h^{*},H^{*}) $ is a balanced pair as in (4). For any
    $(h,H)$,
    \begin{equation}\tP(h, H) \geq \tP(FS(H), H)\geq \tP(FS(H^{*}), H^{*})
    =\tP(h^{*},H^{*}),\end{equation}
    by (17) and Corollary 2 and Lemma 4. We also have 
      $$\tP(h, \Hilb(h)) = \log d + \frac{1}{d} \tL(h), $$
  so taking $H=\Hilb(h)$ in (18) we get
  $$ \tL(h)\geq \tL(h^{*})$$
  as required.
   
   Notice that we also have the symmetrical result to Lemma 4:
   \begin{lem}
     For any $h\in \cK,H\in M$,
     $$  \tP(h,H) \geq \tP(h, \Hilb(h)). $$
     \end{lem}
     This is completely elementary, it just amounts to the standard inequality
     $$   \left(\det Q\right)^{1/d}\leq \frac{{\rm Tr}\ (Q)}{d}, $$
     for $d\times d$ positive Hermitian matrices $Q$. Alternatively, it
     is a simple example of the \lq\lq Kempf-Ness principle''.
   
    Using this Lemma we can obtain a slightly different result: whether
    or not a balanced point exists, the function $\tL$ is bounded below
    on $\cK$ if and only if the function $\tZ$ is bounded below on $M$.
    
    \section{}
    
    We now move on to the second part of the proof of Theorem 1.
    We replace the line bundle $L$ by $L^{k}$ where $k$ is large and positive
     and we let $E_{k}$ be the complete
     linear system $H^{0}(X;L^{k})$. However we work with the same Kahler class
     $\cK$. The whole discussion above goes through, with trivial changes. For
     $h\in \cK$ we define functions $\rho_{k,h}$ on $X$ given by
     $\frac{{\rm dim} E_{k}}{V}\sum \vert s_{\alpha}\vert^{2}$ for a
      $\Hilb(h)$-orthonormal basis of $E_{k}$.

  For a function $f$ on $X$ let us write 
  $$ [f]_{h} = f - \hat{f}, $$
  where $\hat{f}$ is the average value of $f$ defined using the measure
  $d\mu_{h}$. We obtain functionals $\tL_{k}$
   on $\cK$ with
  $$ \delta \tL_{k} = \int_{X} [\Delta' \rho_{k,\phi} + k \rho_{k,\phi}]_{h}
  \ \delta \phi \ d\mu_{h}, $$
  while the definition of the Mabuchi functional is
  $$ \delta \cM = \int_{X} [S]_{h}\ d\mu_{h}. $$
  
  According to the results of Tian, Zelditch, Lu {\it et al.} (see the
  references cited in \cite{kn:don}: also the paper \cite{kn:Ruan} and
  the new paper \cite{kn:new}),
  $$   [\Delta' \rho_{k,h} + k \rho_{k,h}]_{h} \sim \frac{1}{2\pi} k^{n} [S]_{h}, $$
  and the asymptotic
  relation holds uniformly over bounded sets of metrics (in a suitable
  sense), as stated in \cite{kn:don}, Proposition 6.
  By integrating over paths in the space of metrics we deduce
  \begin{prop}
  There are constants $\lambda_{k}\in \bR$ such that
  $$  \frac{2\pi}{k^{n} }\tL_{k} + \lambda_{k} \rightarrow \cM, $$
  as $k\rightarrow \infty$, uniformly over bounded subsets in $\cK$.
  \end{prop}
  Again, note that the Mabuchi functional is only defined up to an additive
  constant anyway.

   Now suppose that there is a constant scalar curvature metric $\omega_{\infty}$
   in the given Kahler class. This corresponds to some metric 
   $h_{\infty}\in
   \cK$. By the main result of \cite{kn:don} there is a sequence
   $h_{k}\in \cK$, for large $k$, realising balanced metrics, and
   converging (in $C^{\infty}$) to $h_{\infty}$ as 
   $k\rightarrow \infty$. Here we use the hypothesis on the automorphism
   group of $(X,L)$, required for the result in \cite{kn:don}. By Theorem
   2, the $h_{k}$ are absolute minima of the functionals $\tL'_{k}=
   \frac{2\pi}{
   k^{n}}
   \tL_{k} + \lambda_{k}$. Now from the fact that the $\cL'_{k}$ converge to
   $\cM$ uniformly over bounded sets of metrics it follows immediately
   that $h_{\infty}$ is an absolute minimum of the Mabuchi functional
   $\cM$.
   
   {\bf Remark.} Notice that Lemma's 4 and 5 imply that the map
   $FS\circ\Hilb:M\rightarrow M$ decreases the functional $\tZ$. It follows
   easily (since $\tZ$ is proper) that iterates of this map converge to the
    minimiser.
   This gives an algorithm for finding balanced metrics numerically
   (when they exist), and hence an algorithm for finding constant scalar
   curvature metrics, when these exist.


\begin{thebibliography}{99}
\bibitem{kn:BM} S. Bando and T. Mabuchi {\em Uniqueness of Einstein-Kahler
metrics modulo connected group actions} In: {\it Algebraic Geometry} Advanced
Studies in Pure Math. 11-40 1987
\bibitem{kn:Chen} X. Chen {\em The space of Kahler metrics} Jour. Differential
Geometry {\bf 56} 189-234 2000
\bibitem{kn:CT} X. Chen and G. Tian {\em Geometry of Kahler metrics and
holomorphic foliation by discs} Preprint DG/0409433
\bibitem{kn:new} X. Dai, K. Liu and X. Ma {\em On the asymptotic expansions
of Bergman kernels} Preprint DG/0404494
\bibitem{kn:D0} S. Donaldson {\em Symmetric spaces, Kahler geometry and
Hamiltonian dynamics} In: {\it Northern California Symplectic Geometry
Seminar} (Eliashberg {\it et al} eds. ) Amer. Math. Soc 1999
\bibitem{kn:don} S.Donaldson {\em Scalar curvature and projective embeddings,
I} Jour. Differential Geometry {\bf 59} 479-522 2001
\bibitem{kn:DK} S. Donaldson and P. Kronheimer {\em The geometry of four-manifolds}
Oxford U.P. 1990
\bibitem{kn:KN} G. Kempf and L. Ness {\em The length of vectors in representation
spaces} In:{\it Algebraic Geometry} Lecture Notes in Math. Springer 
{\bf 732} 233-243 1979
\bibitem{kn:M} T. Mabuchi {\em K-energy maps integrating Futaki invariants}
Tohuku Math. Jour. {\bf 38} 245-257 1986
\bibitem{kn:M0} T. Mabuchi {\em Some symplectic geometry on compact Kahler
manifolds, I} Osaka Jour. Math. {\bf 24} 227-252 1987
\bibitem{kn:M1} T. Mabuchi {\em An energy-theoretic approach to the Hitchin-Kobayashi
correspondence for manifolds, I} To appear in Inventiones Math.
\bibitem{kn:M2} T. Mabuchi {\em An energy-theoretic approach to the Hitchin-Kobayashi
correspondence for manifolds, II} To appear in Jour. Differential Geometry.
\bibitem{kn:M3} T. Mabuchi {\em Stability of extremal Kahler manifolds}
To appear in Osaka Jour. Math.
\bibitem{kn:PS1} D. Phong and J. Sturm {\em Scalar curvature, moment maps
and the Deligns pairing} Preprint DG/02009098
\bibitem{kn:PS2} D. Phong and J. Sturm {\em Stability, energy functionals
and Kahler-Einstein metrics} Preprint DG/0203254
\bibitem{kn:Ruan} W. Ruan {\em Canonical coordinates and Bergman metrics}
Commun. Anal. Geom. {\bf 6} 589-631 1998
\bibitem{kn:Semmes} S. Semmes {\em Complex Monge-Amp\`ere and symplectic
manifolds} Amer. Jour. Math. {\bf 114} 495-550 1992
\bibitem{kn:Z} S. Zhang {\em Heights and reductions of semi-stable varieties}
Compositio Math. {\bf 104} 77-105 1996

\end{thebibliography}
\end{document}